\documentclass[10pt,reqno]{amsart}
\usepackage{amsmath}
\usepackage{amssymb}
\usepackage{amsthm}

\newlength{\defbaselineskip}
\newcommand{\setlinespacing}[1]%
           {\setlength{\baselineskip}{#1 \defbaselineskip}}

\numberwithin{equation}{section}

\newtheorem{thm}{Theorem}[section]

\newtheorem{lem}[thm]{Lemma}

\theoremstyle{definition}

\theoremstyle{remark}
\newtheorem{rem}[thm]{Remark}
\numberwithin{equation}{section}

\begin{document}

\title[Strichartz inequalities for the wave equation]
{Note on Strichartz inequalities for the wave equation with potential}

\author{Seongyeon Kim, Ihyeok Seo and Jihyeon Seok}


\address{Department of Mathematics, Sungkyunkwan University, Suwon 16419, Republic of Korea}
\email{synkim@skku.edu}

\email{ihseo@skku.edu}

\email{jseok@skku.edu}

\begin{abstract}
We obtain Strichartz inequalities for the wave equation with potentials which behave like the inverse square potential $|x|^{-2}$
but might be not a radially symmetric function.
\end{abstract}

\maketitle

\section{Introduction}

Consider the following Cauchy problem for the wave equation with a potential $V(x)$:
\begin{equation}\label{WE}
\begin{cases}
-\partial_t^2 u+\Delta u-V(x)u=0,\\
u(x,0)=f(x), \\
\partial_t u(x,0)=g(x),
\end{cases}
\end{equation}
where $u:\mathbb{R}^{n+1}\rightarrow\mathbb{C}$, $V:\mathbb{R}^{n}\rightarrow\mathbb{C}$
and $\Delta$ is the $n$ dimensional Laplacian.

In this paper we are concerned with the Strichartz inequalities for the wave equation \eqref{WE}.
In the free case $V\equiv0$, the following remarkable estimate was first obtained by Strichartz \cite{S}
in connection with Fourier restriction theory in harmonic analysis:
\begin{equation}\label{oneS}
\|u\|_{L^{\frac{2(n+1)}{n-1}}(\mathbb{R}^{n+1})}
\lesssim \|f\|_{\dot{H}^{1/2}(\mathbb{R}^n)} + \|g\|_{\dot{H}^{-1/2}(\mathbb{R}^n)}, \quad n\ge2,
\end{equation}
where $\dot{H}^\gamma$ denotes the homogeneous Sobolev space equipped with the norm
$\|f\|_{\dot{H}^\gamma} = \| (\sqrt{-\Delta})^\gamma f\|_{L^2}$.
Since then, \eqref{oneS} was extended to mixed norm spaces $L_t^qL_x^r$ as follows (see \cite{LS,KT} and references therein):
\begin{equation*}
\|u\|_{L_t^q(\mathbb{R};L_x^r(\mathbb{R}^n))} \lesssim \|f\|_{\dot{H}^{1/2}(\mathbb{R}^n)} + \|g\|_{\dot{H}^{-1/2}(\mathbb{R}^n)}
\end{equation*}
for $q \ge 2$, $2 \le r < \infty$,
\begin{equation}\label{cond}
\frac{2}{q} + \frac{n-1}{r} \le \frac{n-1}{2}\quad \text{and}\quad\frac{1}{q} + \frac{n}{r} = \frac{n-1}{2}.
\end{equation}
Here the second equality in \eqref{cond} is just the gap condition
and the case $q=r$ recovers the classical result \eqref{oneS}.

Now we turn to the wave equation with a potential.
Several works have treated this potential perturbation of the
free wave equation.
In \cite{BS} the potentials satisfy the decay assumption that $V(x)$ decays like $|x|^{-4-\varepsilon}$ at infinity.
In \cite{C} this assumption is weakened to $|x|^{-3-\varepsilon}$,
which is in turn improved to $|x|^{-2-\varepsilon}$ in \cite{GV}.
In these papers Strichartz type estimates for the corresponding perturbed wave
equations are established.
But the main interest in the equation \eqref{WE} comes from the case where the potential term is homogeneous of degree $-2$
and therefore scales exactly the same as the Laplacian.
For instance, when $V(x)=a|x|^{-2}$ with a real number $a$, the equation \eqref{WE} arises in the study of wave propagation on
conic manifolds \cite{CT}.
We also note that the heat flow for the operator $-\Delta+a|x|^{-2}$ has been studied
in the theory of combustion \cite{VZ}.
In fact, the decay $|V| \sim |x|^{-2}$ was shown to be critical in \cite{GVV} which concerns explicitly the Schr\"odinger case but can be adapted to the wave equation as well. 

For the inverse square potentials $V(x)=a|x|^{-2}$ with $a>-{(n-2)^2}/4$,
Planchon, Stalker and Tahvildar-Zadeh \cite{PST-Z} first obtained the Strichartz inequalities
for the equation \eqref{WE} with radial Cauchy data $f$ and $g$.
Thereafter, this radially symmetric assumption was removed in \cite{BPST-Z}.
More precisely, the range of the admissible exponents $(q,r)$ for the
Strichartz inequalities
\begin{equation*}
\|u\|_{L_t^q(\mathbb{R};\dot{H}^\sigma_r(\mathbb{R}^n))} \lesssim \|f\|_{\dot{H}^{1/2}(\mathbb{R}^n)} + \|g\|_{\dot{H}^{-1/2}(\mathbb{R}^n)}
\end{equation*}
obtained in \cite{PST-Z,BPST-Z} with the gap condition $\sigma=\frac{1}{q} + \frac{n}{r} -\frac{n-1}{2}$ is restricted under
$\frac{2}{q} + \frac{n-1}{r} \le \frac{n-1}{2}$
which is the same as that of the wave equation without potential.
Here,
$\|f\|_{\dot{H}^\sigma_r} = \|(\sqrt{-\Delta})^\sigma f\|_{L^r}$.
In \cite{BPST-Z2} these results were further extended to potentials which behave like the inverse square potential
but might be not a radially symmetric function.
Indeed, the potentials considered in \cite{BPST-Z2} are contained in the weak space, $L^{n/2,\infty}$.

In this paper we consider the Fefferman-Phong class of potentials which is defined for $1\leq p\leq n/2$ by
\begin{equation*}
V \in {\mathcal{F}}^p \quad \Leftrightarrow \quad \|V\|_{\mathcal{F}^p} = \sup_{x \in \mathbb{R}^n, r>0}r^{2-n/p} \left( \int_{B_r(x)} |V(y)|^p dy \right)^{\frac{1}{p}} < \infty,
\end{equation*}
where $B_r(x)$ denotes the ball centered at $x$ with radius $r$.
Note that $L^{n/2}={\mathcal{F}}^{n/2}$ and $a|x|^{-2} \in L^{n/2, \infty} \subsetneq {\mathcal{F}}^p$ if $1 \le p < n/2$.
Our result is the following theorem.

\begin{thm}\label{thm}
Let $n \ge 3$. Let $u$ be a solution to \eqref{WE} with Cauchy data $(f,g)\in\dot{H}^{1/2}\times\dot{H}^{-1/2}$
and potential $V \in {\mathcal{F}^p}$ with small $\|V\|_{\mathcal{F}^p}$ for $p>(n-1)/2$.
Then we have
	\begin{equation}\label{result}
	\|u\|_{L_t^q(\mathbb{R};\dot{H}^\sigma_r(\mathbb{R}^n))}\lesssim
(1+ \|V\|_{\mathcal{F}^p}) \Big(\|f\|_{\dot{H}^\frac{1}{2}} +\|g\|_{\dot{H}^{-\frac{1}{2}}}\Big)
	\end{equation}
for $q > 2$, $2 \le r < \infty$,
\begin{equation}\label{cond2}
\frac{2}{q} + \frac{n-1}{r} \le \frac{n-1}{2}\quad \text{and}\quad\sigma = \frac{1}{q} + \frac{n}{r} -\frac{n-1}{2}.
\end{equation}	
\end{thm}

\begin{rem}
The class of potentials in the theorem is strictly larger than $L^{n/2,\infty}$.
For instance, consider
	\begin{equation*}
	V(x)= f(\frac{x}{|x|}) |x|^{-2}, \quad f \in L^p(S^{n-1}), \quad (n-1)/2< p< n/2,
	\end{equation*}
	which is in ${\mathcal{F}^p}$ but not in $L^{n/2, \infty}$.
\end{rem}

Throughout this paper, the letter $C$ stands for a positive constant which may be different
at each occurrence.
We also denote $A\lesssim B$ to mean $A\leq CB$
with unspecified constants $C>0$.

\


\section{Proof of Theorem \ref{thm}} \label{sec2}

In this section we prove the Strichartz inequalities \eqref{result} by making use of a weighted space-time $L^2$ estimate for the wave equation.

We first consider the potential term as a source term and then write the solution to \eqref{WE} as the sum of the solution to
the free wave equation plus a Duhamel term, as follows:
\begin{equation}\label{sol}
u(x,t) \\= \cos(t\sqrt{-\Delta})f + \frac{\sin(t\sqrt{-\Delta})}{\sqrt{-\Delta}}g +
\int_{0}^{t} \frac{\sin((t-s)\sqrt{-\Delta})}{\sqrt{-\Delta}}\big(V(\cdot)u(\cdot, s)\big) ds.
\end{equation}
By the classical Strichartz inequalities for the wave equation (see e.g. \cite{KT}), we see
\begin{equation}\label{homo}
\| e^{it\sqrt{-\Delta}}f\|_{L^q_t \dot{H}^\sigma_r} = \| (\sqrt{-\Delta})^\sigma e^{it\sqrt{-\Delta}}f\|_{L^q_t L^r_x} \lesssim \|f\|_{\dot{H}^{\frac{1}{2}}}
\end{equation}
for $(q,r)$ satisfying $q \ge 2$, $2 \le r < \infty$ and the condition \eqref{cond2}.
Applying \eqref{homo} to \eqref{sol}, we get
\begin{align*}
\|u\|_{{L^q_t \dot{H}^\sigma_r}} &
\lesssim \|f\|_{\dot{H}^{\frac{1}{2}}} + \|g\|_{\dot{H}^{-\frac{1}{2}}} + \bigg\|\int_{0}^{t} \frac{\sin((t-s)\sqrt{-\Delta})}{\sqrt{-\Delta}}\big(V(\cdot)u(\cdot, s)\big) ds\bigg\|_{L^q_t \dot{H}^\sigma_r}
\end{align*}
for the same $(q,r)$.

Now it remains to show that
\begin{equation*}
\bigg\|\int_{0}^{t} \frac{\sin((t-s)\sqrt{-\Delta})}{\sqrt{-\Delta}}\big(V(\cdot)u(\cdot, s)\big) ds\bigg\|_{L^q_t \dot{H}^\sigma_r} \lesssim  \|V\|_{\mathcal{F}^p}  \big(\|f\|_{\dot{H}^{\frac{1}{2}}} + \|g\|_{\dot{H}^{-\frac{1}{2}}}\big)
\end{equation*}
for $(q,r)$ satisfying $q> 2$, $2 \le r < \infty$ and the condition \eqref{cond2}.
By duality, it is sufficient to show that
\begin{equation}\label{d-inhom}
\begin{split}
\bigg< (\sqrt{-\Delta})^\sigma \int_{0}^{t} \frac{\sin((t-s)\sqrt{-\Delta})}{\sqrt{-\Delta}}&\big(V(\cdot)u(\cdot, s)\big) ds , G \bigg>_{x,t} \\
& \lesssim \|V\|_{\mathcal{F}^p} \big(\|f\|_{\dot{H}^{\frac{1}{2}}} + \|g\|_{\dot{H}^{-\frac{1}{2}}}\big){\|G\|}_{L^{q^\prime }_t L^{r^\prime }_x}.
\end{split}
\end{equation}
The left-hand side of \eqref{d-inhom} is equivalent to
\begin{align}\label{eq1}
\qquad \qquad \int_{\mathbb{R}}\int_{0}^{t}&{\big< (\sqrt{-\Delta})^{\sigma-1}\sin((t-s)\sqrt{-\Delta})\big(V(\cdot)u(\cdot, s)\big), G\, \big>}_{x} dsdt \nonumber
\\&= \int_{\mathbb{R}}\int_{0}^{t}{\big<  Vu, (\sqrt{-\Delta})^{\sigma-1} \sin((t-s)\sqrt{-\Delta})G\, \big>}_{x} dsdt \nonumber
\\&=\bigg<  V^{1/2} u, V^{1/2} (\sqrt{-\Delta})^{\sigma-1} \int_{s}^{\infty} \sin((t-s)\sqrt{-\Delta})G\,dt \bigg>_{x, s}.
\end{align}
Using H\"older's inequality, \eqref{eq1} is bounded by
\begin{equation*}
\|u\|_{L^2_{x,s}(|V|)}  \Big\|(\sqrt{-\Delta})^{\sigma-1}\int_{s}^{\infty}\sin((t-s)\sqrt{-\Delta})\,G\,dt \Big\|_{L^2_{x,s}(|V|)}.
\end{equation*}
Here, $L^2(|V|)$ denotes a weighted space equipped with the norm
$$\|h\|_{L_{x,t}^2(|V|)}= \bigg(\int_{\mathbb{R}^{n+1}} |h(x,t)|^2 |V(x)| dxdt\bigg)^{\frac{1}{2}}.$$
We will show that
\begin{equation}\label{d1}
\|u\|_{L^2_{x,t}(|V|)} \lesssim \|V\|^{1/2}_{\mathcal{F}^p}\big(\|f\|_{\dot{H}^\frac{1}{2}} + \|g\|_{\dot{H}^{-\frac{1}{2}}}\big)
\end{equation}
and
\begin{equation}\label{d2}
\Big\|(\sqrt{-\Delta})^{\sigma-1}\int_{t}^{\infty}\sin((t-s)\sqrt{-\Delta})\,G\,ds \, \Big\|_{L^2_{x,t}(|V|)} \lesssim \|V\|^{1/2}_{\mathcal{F}^p}{\|G\|}_{L^{q^\prime }_t L^{r^\prime }_x}
\end{equation}
for $(q,r)$ satisfying the same conditions in the theorem.
Then the desired estimate \eqref{d-inhom} is proved.

To show \eqref{d1}, we use the following lemma which is a particular case of
Proposition 2.3 and 4.2 in \cite{RV}.

\begin{lem}\label{lem1}
Let $n \ge 3$.
Assume that $V\in\mathcal{F}^p$ for $p > (n-1)/2$.
Then we have
	\begin{equation}\label{cos}
	\|\cos (t\sqrt{-\Delta}) f\|_{L_{x,t}^2(|V|)} \lesssim \|V\|^{1/2}_{\mathcal{F}^p} \|(\sqrt{-\Delta})^{1/2}f\|_{L^2},
	\end{equation}
	\begin{equation}\label{sin}
	\bigg\|\frac{\sin t\sqrt{-\Delta}}{\sqrt{-\Delta}} g\bigg\|_{L_{x,t}^2(|V|)}
 \lesssim \|V\|^{1/2}_{\mathcal{F}^p} \|(\sqrt{-\Delta})^{-1/2}g\|_{L^2},
	\end{equation}	
and
	\begin{equation*}
	\bigg\|\int_{0}^{t} \frac{\sin((t-s)\sqrt{-\Delta})}{\sqrt{-\Delta}}F(\cdot, s)\, ds\bigg\|_{L_{x,t}^2(|V|)}
 \lesssim \|V\|_{\mathcal{F}^p} \|F\|_{L_{x,t}^2(|V|^{-1})}.
	\end{equation*}	
\end{lem}

Indeed, applying Lemma \ref{lem1} to \eqref{sol}, we see
\begin{equation}\label{eq2}
\begin{split}
\|u\|_{L^2_{x,t}(|V|)} \lesssim \|V\|^{1/2}_{\mathcal{F}^p} \big(\|f\|_{\dot{H}^\frac{1}{2}} + \|g\|_{\dot{H}^{-\frac{1}{2}}}\big)+ \|V\|_{\mathcal{F}^p}\|u\|_{L^2_{x,t}(|V|)}.
\end{split}
\end{equation}
Since we are assuming that $\|V\|_{\mathcal{F}^p}$ is small enough, the last term on the right-hand side of \eqref{eq2}
can be absorbed into the left-hand side.
Hence, we get the first estimate \eqref{d1}.
To obtain the second estimate \eqref{d2}, we first note that the first two estimates \eqref{cos} and \eqref{sin} in Lemma \ref{lem1}
directly imply
\begin{equation}\label{lem1-result}
\big\|(\sqrt{-\Delta})^\gamma e^{it\sqrt{-\Delta}} f\big\|_{L_{x,t}^2(|V|)}
\lesssim \|V\|^{1/2}_{\mathcal{F}^p} \big\|(\sqrt{-\Delta})^{\gamma+1/2}f\big\|_{L^2}
\end{equation}
for $\gamma\in\mathbb{R}$.
Using \eqref{lem1-result} and the dual estimate of \eqref{homo}, we then have
\begin{align*}
\Big\|(\sqrt{-\Delta})^{\sigma-1}\int_{\mathbb{R}}\sin((t-s)\sqrt{-\Delta})&\,G\,ds \Big\|_{L^2_{x,t}(|V|)}\\
\nonumber&\lesssim \Big\|({\sqrt{-\Delta}})^{\sigma-1}\int_{\mathbb{R}}e^{i(t-s)\sqrt{-\Delta}} \,G\, ds\Big\|_{L^2_{x,t}(|V|)} \nonumber\\&\lesssim \|V\|^{1/2}_{\mathcal{F}^p} \Big\|({\sqrt{-\Delta}})^{\sigma-1/2}\int_{\mathbb{R}} e^{-is\sqrt{-\Delta}}\,G\, ds\Big\|_{L^2} \nonumber\\&\lesssim \|V\|^{1/2}_{\mathcal{F}^p}\|G\|_{L^{q^\prime }_tL^{r^\prime }_x}
\end{align*}
for $(q,r)$ satisfying $q \ge 2$, $2 \le r < \infty$ and the condition \eqref{cond2}.
Here we are going to use the following Christ-Kiselev lemma (\cite{CK}) to conclude that
\begin{equation}\label{ck}
\Big\|(\sqrt{-\Delta})^{\sigma-1}\int_{-\infty}^{t}\sin((t-s)\sqrt{-\Delta})\, G\,ds \Big\|_{L^2_{x,t}(|V|)} \lesssim ||V||^{1/2}_{\mathcal{F}^p}||G||_{L^{q^\prime }_tL^{r^\prime }_x}
\end{equation}
if $2>q'$.

\begin{lem}
Let $X$ and $Y$ be two Banach spaces and let $T$ be a bounded linear operator from $L^\alpha(\mathbb{R};X)$ to $L^\beta(\mathbb{R};Y)$
such that
$$Tf(t)=\int_{\mathbb{R}} K(t,s)f(s)ds.$$
Then the operator
$$\widetilde{T}f(t)=\int_{-\infty}^t K(t,s)f(s)ds$$
has the same boundedness when $\beta>\alpha$, and $\|\widetilde{T}\|\lesssim\|T\|$.
\end{lem}

The desired estimate \eqref{d2} follows directly from \eqref{ck} by changing some variables.
This completes the proof.

 \ 
 
\noindent {\bf{Acknowledgements.}}
This research was supported by NRF-2019R1F1A1061316. 
The authors would like to thank Y. Koh for discussions on related issues.


\begin{thebibliography}{99}

\bibitem{BS} M. Beals and W. Strauss, \textit{$L^p$ estimates for the wave equation with a
potential},
Comm. Partial Differential Equations 18(1993), 1365-1397.
	
\bibitem{BPST-Z} N. Burq, F. Planchon, J. G. Stalker and A. S. Tahvildar-Zadeh, \textit{Strichartz estimates for the wave and Schr\"odinger equations with the inverse-square potential}, J. Funct. Anal. 203(2003), 519-549.

\bibitem{BPST-Z2} N. Burq, F. Planchon, J. G. Stalker and A. S. Tahvildar-Zadeh, \textit{Strichartz estimates for the wave and Schr\"odinger equations with potentials of critical decay}, Indiana Univ. Math. J. 53(2004), 1665-1680.

\bibitem{CT} J. Cheeger and M. Taylor, \textit{On the diffraction of waves by conical singularities. I},
Comm. Pure Appl. Math. 35(1982), 275-331.

\bibitem{CK} M. Christ and A. Kiselev, \textit{Maximal operators associated to filtrations}, J. Funct. Anal. 179(2001), 409-425.

\bibitem{C} S. Cuccagna, \textit{On the wave equation with a potential},
Comm. Partial Differential Equations 25(2000), 1549-1565.

\bibitem{GV} V. Georgiev and N. Visciglia, \textit{Decay estimates for the wave equation with potential},
Comm. Partial Differential Equations 29(2003), 101-153.

\bibitem{GVV} M. Goldberg, L. Vega and N. Visciglia, \textit{Counterexamples of Strichartz inequalities for Schr\"odinger equations with repulsive potentials}, Int. Math. Res. Not. 2006, Art. ID 13927, 16pp.

\bibitem{KT} M. Keel and T. Tao, \textit{Endpoint strichartz estimates}, Amer. J. math. 120(1998), 955-980.


\bibitem{LS} H. Lindblad and C. D. Sogge, \textit{On existence and scattering with minimal regularity for semilinear wave equations}, J. Funct. Anal. 130(1995), 357-426.

\bibitem{PST-Z} F. Planchon, J. G. Stalker and A. S. Tahvildar-Zadeh, \textit{$L^p$ estimates for the wave equation with the inverse-square potential}, Discrete Contin. Dyn. Syst. 9(2003), 427-442.

\bibitem{RV} A. Ruiz and L. Vega, \textit{Local regularity of solutions to wave equation with time-dependent potentials}, Duke Math. J. 76(1994), 913-940.

\bibitem{S} R. Strichartz, \textit{Restrictions of Fourier transforms to quadratic surfaces and decay of solutions of wave equations},
Duke Math. J. 44(1977), 705-714.

\bibitem{VZ} J. L. Vazquez and E. Zuazua, \textit{The Hardy inequality and the asymptotic behaviour
of the heat equation with an inverse-square potential},
J. Funct. Anal. 173(2000), 103-153.

\end{thebibliography}
\end{document}